\documentclass[a4paper,oneside,twocolumn,english]{amsart}
\usepackage[T1]{fontenc}
\usepackage[latin9]{inputenc}
\usepackage{xcolor}
\usepackage{babel}
\usepackage{array}
\usepackage{varioref}
\usepackage{float}
\usepackage{multirow}
\usepackage{amstext}
\usepackage{amsthm}
\usepackage{amssymb}
\usepackage{graphicx}

\makeatletter

\pdfpageheight\paperheight
\pdfpagewidth\paperwidth

\providecommand{\tabularnewline}{\\}
\floatstyle{ruled}
\newfloat{algorithm}{tbp}{loa}
\providecommand{\algorithmname}{Algorithm}
\floatname{algorithm}{\protect\algorithmname}

\numberwithin{equation}{section}
\numberwithin{figure}{section}
\theoremstyle{plain}
\newtheorem{thm}{\protect\theoremname}
\theoremstyle{remark}
\newtheorem{rem}[thm]{\protect\remarkname}

\usepackage{geometry}
\usepackage{url}

\usepackage{mathtools,amssymb,lipsum}

\usepackage{cuted}
\usepackage{hyperref}
\usepackage{amsthm}

\makeatother

\providecommand{\remarkname}{Remark}
\providecommand{\theoremname}{Theorem}

\begin{document}
\title{Velocity-aided IMU-based Attitude Estimation}
\author{Mehdi Benallegue, Abdelaziz Benallegue, Rafael Cisneros, Yacine Chitour }
\begin{abstract}
This paper addresses the problem of estimating the attitude of a rigid
body, which is subject to high accelerations and equipped with inertial
measurement unit (IMU) and sensors providing the body velocity (expressed
in the reference frame attached to the body). That issue can be treated
differently depending on the level of confidence in the measurements
of the magnetometer of the IMU, particularly with regard to the observation
of the inclination component with respect to the vertical direction,
rendering possible to describe the interaction with gravity. Two cases
are then studied: either (i) the magnetometer is absent and only the
inclination can be estimated, (ii) the magnetometer is present, giving
redundancy and full attitude observability. In the latter case, the
presented observer allows to tune how much the inclination estimation
is influenced by the magnetometer. All state estimators are proposed
with proof of almost global asymptotic stability and local exponential
convergence. Finally, these estimators are compared with state-of-the-art
solutions in clean and noisy simulations, allowing recommended solutions
to be drawn for each case.\thanks{M. Benallegue and R. Cisneros are with Humanoid Research Group, National
Institute of Advanced Industrial Science and Technology (AIST), Tsukuba,
Ibaraki, Japan. A. Benallegue is with Université Paris-Saclay, UVSQ,
Laboratoire d'Ingénierie des Systèmes de Versailles, 78124, Vélizy
-Villacoublay, France and JRL-AIST (Joint Robotics Laboratory), Tsukuba,
Ibaraki, Japan. Y. Chitour is with Université Paris-Saclay, CNRS,
CentraleSupélec,  Laboratoire des Signaux et Systèmes,
91190, Gif-sur-Yvette, France. \texttt{\small{}mehdi.benallegue@aist.go.jp,
abdelaziz.benallegue@uvsq.fr, rafael.cisneros@aist.go.jp, yacine.chitour@l2s.centralesupelec.fr}}
\end{abstract}

\maketitle

\section{Introduction}

The orientation, or attitude, of a mechanical system in the world
is often an important part of its dynamical state and constitutes
sometimes the most important variable determining the dynamics such
as in the case of drones~\cite{franklin1994feedback} or legged robots~\cite{Wieber2016}.
While it is possible for fixed-based robots to reconstruct the orientation
of any link using the joint position, this is not possible for mobile
robots, and specifically floating-base ones. However, a set of sensors
is usually dedicated to the estimation of the orientation. These sensors,
usually grouped in a set called inertial measurement units, measure
the linear acceleration, including the gravitational one, the angular
velocity and sometimes the magnetic field measurement, all expressed
in the frame of the sensor. Thanks to this set of measurements an
efficient estimation for the attitude can be built, but only when
the system has negligible linear accelerations compared to gravity
\cite{mahony2008nonlinear,martin2008invariant}.

The constraint that the inertial accelerations must be negligible
compared to gravity can be limiting to low dynamics motion or simply
impossible to hold, especially when the system is subject to impacts
such as during the case of bipedal walking. A dynamical model of the
system can be used to predict the accelerations and compensate for
them. This prediction can be based on the forces models, either in
the case of unmanned aerial vehicles~\cite{Martin2010ICRA,Martin2016cdc}
or legged robots~\cite{Bloesch-RSS-12,mifsud:hal-01142399}. However
this solution is specific to every dynamical system and requires to
identify many dynamical parameters.

Another solution is to ``aid'' the inertial measurement unit (IMU)
with independent measurements able to clear the acceleration ambiguity,
such as the position in the world frame provided by GPS \cite{Hua2010cep,Roberts2011cdc}
or linear velocity, either expressed in the world frame~\cite{Martin2008ifac},
the local frame of the sensor~\cite{Hua2016automatica,allibert2016velocity,Benallegue2017humanoids},
or a mixture of both~\cite{Hua2017cdc}, sometimes while reconstructing
the velocity itself~\cite{allibert2014estimating}. For instance,
the velocity-aided solution allows to reconstruct the attitude efficiently
with proven Lyapunov convergence. In this work, we similarly assume
that the velocity in the sensor local frame is available. This velocity
can be provided by a sensor such as Doppler effect radars. It can
be provided also by the measurements of the gyrometer in the presence
of a known anchor in the environment. This is for example the case
of humanoid robots in contact with the environment, because the contact
point position and velocity in the sensor frame are known~\cite{Benallegue2017humanoids}.
There are a few estimators considering the same case. In~\cite{allibert2016velocity}
a velocity-aided estimator with proof of convergence has been presented,
but there were possible cases of singularities if the scaling factor
reaches zero. In \cite{Hua2016automatica} two estimators with proof
of convergence have been presented, including one without gain condition,
but we presented in \cite{Benallegue2017humanoids} a slight improvement
of the estimator proposed in \cite{Hua2016automatica} where better
performance and simple convergence analysis are obtained. In~\cite{Martin2016arxiv},
a global estimator has been presented, but the globality has been
reached at the cost of breaking the normality constraint of the gravity
and magnetic field direction vectors. It presented also a projected
rotation matrix on SO(3), which can be discontinuous in case of singularities.

We propose here to extend the work presented in \cite{Benallegue2017humanoids}
and at the same time the work of \cite{Hua2016automatica} and \cite{Martin2016arxiv}
with a set of improved estimators. The choice of the estimator should
be made mainly according to the availability and quality of the magnetometer
measurements. The estimators include (i) an estimator for the tilt
(roll and pitch angles) that does not use the magnetometer and (ii)
a set of complete attitude estimators based on our confidence in the
magnetometer. We first introduce a new tilt observer called a \textquotedbl two-step
state observer\textquotedbl{} which operates in two steps: the first
one provides an intermediate estimate $\hat{x}_{2}'\in\mathbb{R}^{3}$
of $x_{2}$ while the second step furnishes the recommended estimate
$\hat{x}_{2}\in\mathbb{S}^{2}$ of $x_{2}$ based on $\hat{x_{2}}'$.
The expected efficiency of this estimator is that its two constitutive
steps are decoupled: the first one insures global exponential convergence
of $\hat{x}_{2}'$ towards $x_{2}$ while the second one is an $\mathbb{S}^{2}$-constrained
complementary-filter estimation also providing better robustness to
disturbances. Finally, this estimator is used to obtain the complete
rotation in the same way as the authors quoted. All full rotation
estimators share a common simple structure, proof of asymptotic convergence
and good overall performance. The quality of the estimate is evaluated
by the comparative simulation with a state-of-the-art solution.

\section{Problem statement\label{sec:Problem-statement}}

\subsection{Frames and measurements}

The problem we address is the estimation of the attitude of accelerated
rigid body vehicles moving in 3D-space. We denote $\mathcal{W}$ the
world frame and $\mathcal{L}$ the local frame of the sensor. This
attitude estimation has to rely on an IMU consisting in a three-axial
accelerometer, gyrometer and magnetometer, and using a measurement
of the velocity of the sensor. The accelerometer provides $y_{a}$
the sum of the gravitational field and the linear acceleration of
the sensor, the gyrometer provides $y_{g}$ measuring the angular
velocity $\omega$ of the IMU and the magnetometer provides $y_{m}$
the measurement of the unit vector $m$ along the Earth's magnetic
field, all these signals are expressed in the sensor frame $\mathcal{L}$.
The velocity sensor provides $y_{v}$ the linear velocity $v$ of
the local frame $\mathcal{L}$ with respect to the world $\mathcal{W}$,
but expressed in $\mathcal{L}$ 
\begin{align}
y_{v}= & v,\\
y_{g}= & \omega,\label{eq:mesure_IMU1}\\
y_{a}= & S(\omega)v+\dot{v}+g_{0}R^{T}e_{z},\label{eq:mesure_IMU2}\\
y_{m}= & R^{T}m,\label{eq:mesure_IMU3}
\end{align}
$R$, $g_{0}$, $e_{z}$ and $m$ are respectively the orientation
of the IMU with respect to the world, the standard gravity constant,
a unit vector collinear with the gravitational field, expressed in
$\mathcal{W}$ and directed upward, and a unit vector aligned with
the earth's magnetic field expressed in $\mathcal{W}$. Finally $\omega$
is the angular velocity of the sensor expressed in $\mathcal{L}$
such that 
\begin{equation}
\dot{R}=RS(\omega),\label{eq:kinematics}
\end{equation}
where the orientation $R$ is the attitude we wish to estimate using
these available measurements and the function $S$ is the skew-symmetric
matrix operator allowing to perform cross-product.

\subsection{State definition and dynamics}

Let us define the following state variables: 
\begin{eqnarray}
x_{1} & \overset{\Delta}{=} & v,\label{eq:x1}\\
x_{2} & \overset{\Delta}{=} & R^{T}e_{z},\label{eq:x2}\\
x_{3} & \overset{\Delta}{=} & R^{T}m,\label{eq:R_c}
\end{eqnarray}
where $x_{1}\in\mathbb{R}^{3}$, $x_{2}\in\mathbb{S}^{2}$ and $x_{3}\in\mathbb{S}^{2}$,
with the set $\mathbb{S}^{2}\subset\mathbb{R}^{3}$ being the unit
sphere centered at the origin, and defined as 
\begin{equation}
\mathbb{S}^{2}\overset{\Delta}{=}\left\{ x\in\mathbb{R}^{3}/\left\Vert x\right\Vert =1\right\} .
\end{equation}

The variables $x_{1}$ and $x_{3}$ are measured using $y_{v}$ and
$y_{m}$, even if they are noisy. On the contrary, $x_{2}$ is the
tilt which cannot be obtained algebraically from the measurements.

From equations (\ref{eq:mesure_IMU2}) and (\ref{eq:x1}) we get 
\begin{align}
\dot{x}_{1}= & -S(\omega)x_{1}+y_{a}-g_{0}R^{T}e_{z}.\label{eq:x1-dynamics}
\end{align}

This, together with the time-differentiation of $x_{2}$ and $x_{3}$
using equation (\ref{eq:kinematics}), provide us with the following
state dynamic equations 
\begin{equation}
\begin{cases}
\dot{x}_{1} & =-S(\omega)x_{1}+y_{a}-g_{0}x_{2},\\
\dot{x}_{2} & =-S(\omega)x_{2},\\
\dot{x}_{3} & =-S(\omega)x_{3}.
\end{cases}\label{eq:dynamics}
\end{equation}

The system (\ref{eq:dynamics}) is suitable for the observer synthesis.

In the next sections we show how to reconstruct the attitude in two
cases: 
\begin{itemize}
\item In Section~\ref{sec:Tilt-estimation}: There is no magnetometer and
only the tilt defined by $x_{2}=R^{T}e_{z}$ can be estimated. 
\item In Section~\ref{sec:Invariant}: The magnetometer is available and
the full rotation matrix $R$ can be reconstructed. The presence of
the magnetometer provides us with some redundancy but the magnetometer's
measurements may be not reliable. This requires us to be able to tune
how much the magnetometer interferes with tilt estimation.
\end{itemize}

\subsection{Basic facts and notation}

We next introduce notations and recall basic properties used in the
developments below where $v$, $w$ and $u$ are vectors and $R\in SO(3)$
a rotation matrix 
\begin{align}
S\left(v\right)S\left(w\right) & =wv^{T}-\left(v^{T}w\right)I,\label{eq:propAH}\\
S\left(v\right)S\left(w\right)S(v) & =-\left(v^{T}w\right)S(v),\label{eq:propAM}\\
R\left(\prod_{i=1}^{n}S\left(v_{i}\right)\right)R^{T} & =\prod_{i=1}^{n}S\left(Rv_{i}\right),\label{eq:propB}\\
S\left(S\left(v\right)w\right) & =S\left(v\right)S\left(w\right)-S\left(w\right)S\left(v\right),\label{eq:propC}\\
 & =wv^{T}-vw^{T},\label{eq:propCM}\\
S^{3}\left(v\right) & =-\Vert v\Vert^{2}S\left(v\right),\label{eq:propE}
\end{align}
where $I$ denotes the $3\times3$ identity matrix.

For $n\geq1$, define $\varUpsilon_{n}\overset{\Delta}{=}\mathbb{R}^{3n}\times\mathbb{S}_{e_{z}}$
and $\varUpsilon_{n}^{*}\overset{\Delta}{=}\mathbb{R}^{3n}\times\mathbb{S}_{e_{z}}^{*}$
with $\mathbb{S}_{e_{z}}\overset{\Delta}{=}\left\{ z\in\mathbb{R}^{3}|\left(e_{z}-z\right)\in\mathbb{S}^{2}\right\} $
and $\mathbb{S}_{e_{z}}^{*}\overset{\Delta}{=}\mathbb{S}_{e_{z}}\setminus\{2e_{z}\}$. 

Finally we will say that a dynamical system $(D)\ \dot{x}=f(x)$ defined
on a differential manifold $X$ is \emph{almost globally asymptotically
stable with respect to an equilibrium point $x_{0}$ of $f$}{} \textit{if
$(D)$ is (Lyapunov) locally stable with respect to $x_{0}$ and there
exists an open an dense subset $X_{1}$ of $X$ such that every trajectory
starting in $X_{1}$ converges asymptotically to $x_{0}$.}

\section{Tilt estimation\label{sec:Tilt-estimation}}

In this section we ignore the signals of the magnetometer, either
because it is unavailable or because the magnetic field is not steady.
In this case the whole orientation cannot be observed. Nevertheless,
we show hereinafter how we can have an efficient estimation of tilt
$x_{2}=R^{T}e_{z}$.

We present a novel observer called ``Two-steps state observer''
which is designed in two steps: the first step provides $\hat{x}'_{2}\in\mathbb{R}^{3}$
an intermediate estimate of $x_{2}$; and the second provides $\hat{x}_{2}\in\mathbb{S}^{2}$
the recommended estimate of $x_{2}$ based on the intermediate one
$\hat{x}'_{2}$. The expected efficiency of this estimator is that
it relies on two stages, the first independent one is given by $\hat{x}_{2}^{\prime}$
which is globally exponentially converging to $x_{2}$ in an efficient
way, and the second is a $\mathbb{S}^{2}$-constrained complementary-filter
estimation providing continuity and better robustness to disturbances.
Indeed, the global exponential convergence of the error of the first
stage actually leads to the violation of the normality constraint
of $R^{T}e_{z}$. Furthermore, the simple normalization of this estimation
may lead undefined output and unbounded time-derivatives when the
norm is close to zero. This can cause a problem when continuity of
the estimation is required. A simple solution is to add the second
stage to maintain the constraint of the tilt estimation in $\mathbb{S}^{2}$
while keeping bounded velocities.

\subsection{Two-steps first order state observer designed in \texorpdfstring{{}$\mathbb{R}^{3}\times\mathbb{S}^{2}$
}{R3xS2}}

\label{ssection-TSFO}

The simplest two-steps estimator can be described as follows, 
\begin{equation}
\begin{cases}
\dot{\hat{x}}_{1} & =-S(y_{g})\hat{x}_{1}+y_{a}-g_{0}\hat{x}_{2}',\\
\hat{x}_{2}^{\prime} & =-\frac{\alpha_{1}}{g_{0}}\left(y_{v}-\hat{x}_{1}\right),\\
\dot{\hat{x}}_{2} & =-S\left(y_{g}-\gamma S(\hat{x}_{2})\hat{x}_{2}^{\prime}\right)\hat{x}_{2},
\end{cases}\label{eq:observer-1st-order}
\end{equation}
where $\alpha_{1}$ and $\gamma$ are positive scalar gains.

If the initial value of $\hat{x}_{2}$ is in $\mathbb{S}^{2}$, then
the dynamics of the last equation ensures that the norm of this vector
remains in time constantly equal one. The initial value for $\hat{x}_{1}$
on the other side could be anywhere in $\mathbb{R}^{3}$.

Using the estimation errors defined as $\tilde{x}_{2}^{\prime}\overset{\Delta}{=}x_{2}-\hat{x}_{2}^{\prime}=x_{2}+\frac{\alpha_{1}}{g_{0}}\left(x_{1}-\hat{x}_{1}\right)=p_{1}$
and $\tilde{x}_{2}\overset{\Delta}{=}x_{2}-\hat{x}_{2}$, and equation
(\ref{eq:dynamics}) we get the error dynamics as 
\begin{equation}
\begin{cases}
\dot{p}_{1} & =-S(\omega)p_{1}-\alpha_{1}p_{1},\\
\dot{\tilde{x}}_{2} & =-S(\omega)\tilde{x}_{2}+\gamma S^{2}\left(\hat{x}_{2}\right)\left(\tilde{x}_{2}-p_{1}\right).
\end{cases}\label{eq:error_dynamics-1st-order}
\end{equation}

To run the analysis of errors, we set $z_{p_{1}}=Rp_{1}$ and $z_{2}=R\tilde{x}_{2}$.
Noticing $R\hat{x}_{2}=e_{z}-z_{2}$, one gets 
\begin{equation}
\begin{cases}
\dot{z}_{p_{1}} & =-\alpha_{1}z_{p_{1}},\\
\dot{z}_{2} & =\gamma S^{2}\left(e_{z}-z_{2}\right)\left(z_{2}-z_{p_{1}}\right).
\end{cases}\label{eq:error_dynamics-2-1st-order}
\end{equation}

This new error dynamics is autonomous and defines a time-invariant
ordinary differential equation (ODE). If one defines the state $\xi_{1}\overset{\Delta}{=}\left(z_{p_{1}},z_{2}\right)\in\varUpsilon_{1}$
one can write (\ref{eq:error_dynamics-2-1st-order}) as $\dot{\xi}_{1}=F_{1}\left(\xi_{1}\right)$
where $F_{1}$ gathers the right-hand side of (\ref{eq:error_dynamics-2-1st-order})
and defines a smooth vector field on $\varUpsilon_{1}$.

We now turn to the convergence analysis of (\ref{eq:error_dynamics-2-1st-order})
and we get the following. 
\begin{thm}
\label{th1}The time-invariant ODE defined by (\ref{eq:error_dynamics-2-1st-order})
verifies the following 
\end{thm}

\begin{enumerate}
\item The state space is equal to $\varUpsilon_{1}$, it admits two equilibrium
points namely the origin $(0,0)$ and $(0,2e_{z})$ and all trajectories
of (\ref{eq:error_dynamics-2-1st-order}) converge to one of the two
equilibrium points. 
\begin{enumerate}
\item \emph{The system }(\ref{eq:error_dynamics-2-1st-order})\emph{ is
almost globally asymptotically stable with respect to the origin,
which is locally exponentially stable. } 
\item For every compact set $K$ of $\varUpsilon_{1}^{*}$ and positive
number $\varrho>0$, there exists $(\alpha_{1},\gamma)$ such that
trajectories of (\ref{eq:error_dynamics-2-1st-order}) starting in
$K$ converge exponentially to the origin with an exponential rate
larger than or equal to $\varrho$.  
\end{enumerate}
\end{enumerate}
The proof of the theorem is given in Section~\ref{app:proof-th1}.
\begin{rem}
\label{rem1} The estimator for the tilt $x_{2}$ operates in two
\emph{decoupled} steps: the first one shows that the artificial state
$\hat{x}_{2}^{\prime}$ estimates $x_{2}$ (the dynamics of the error
term ${z}_{p_{1}}$ is independent of the rest of the system dynamics)
and then, in the second step, one brings back $\hat{x}_{2}^{\prime}$
on $\mathbb{S}^{2}$ through $\hat{x}_{2}$. 
\end{rem}

We must now compare our tilt estimator with previous works. We start
by considering the seminal work \cite{Hua2016automatica} and the
tilt estimator introduced there, that we actually recall below in
(\ref{eq:RxS_observer_Hua}). The price to pay in the present paper
with respect to (\ref{eq:RxS_observer_Hua}) is the extra state $\hat{x}_{2}^{\prime}$
but it has the advantage of not being constrained to $\mathbb{S}^{2}$
anymore. This is clearly put forward when one compares the error dynamics
given by (\ref{eq:error_dynamics-2-1st-order}) and (\ref{eq:RxS_error_dynamics_Hua}).
The decoupling in (\ref{eq:error_dynamics-2-1st-order}) between the
errors $z_{p_{1}}$ and $z_{2}$ not only allows one to have better
convergence results with respect to (\ref{eq:RxS_observer_Hua}) (much
simpler convergence analysis, no conditions on the gains and arbitrary
rate of exponential convergence) but also to improve the robustness
of the estimation to noise.

The second major reference for tilt estimation is that of \cite{Martin2016arxiv}
where the authors provide an estimator in $\mathbb{R}^{3}\times\mathbb{R}^{3}$
instead of $\mathbb{R}^{3}\times\mathbb{S}^{2}$ (cf. the variables
$\hat{v}$ and $\hat{\gamma}$). The error system turns out to be
(essentially) linear and time invariant with, therefore, the best
convergence properties. However, the tilt estimator $\hat{\gamma}$
does not belong to $\mathbb{S}^{2}$ and that may create singularity
issues (i.e., $\hat{\gamma}$ may be equal to zero during the estimation
or becomes collinear to another important unit vector) when one is
interested by tilt estimation only or uses such an estimator to get
an estimator of the total rotation $R$ by, for instance, the TRIAD
method~\cite{shuster1981jgc}.

In conclusion, our tilt estimator $\hat{x}_{2}$ combines the good
convergence properties of the tilt estimator of \cite{Martin2016arxiv}
with the fact that it remains on $\mathbb{S}^{2}$, like the tilt
estimator of \cite{Hua2016automatica}. 

\subsection{Two-steps \texorpdfstring{$n^{th}$}{n-th} order state observer
designed in \texorpdfstring{{}{}{$\mathbb{R}^{3n}\times\mathbb{S}^{2}$}}{R3nxS2} }

\label{subsec:State-observer-in-R3xR3xS2}

An interesting way to comprehend the estimator (\ref{eq:observer-1st-order})
is by noting that the dynamics of $\hat{x}_{2}^{\prime}$ has a first
order exponential convergence to $x_{2}$ and that the dynamics of
$\hat{x}_{2}$ is a complementary filter of $\hat{x}_{2}^{\prime}$.
Therefore, we can extend this feature to higher order of exponential
convergence while keeping the same two-steps structure tilt estimator.
Let $n\geq2$ be an integer. The $n$-th order observer is designed
on $\mathbb{R}^{3n}\times\mathbb{S}^{2}$, where $n$ is the order
of the filter for the first step of the estimator. Increasing the
order of linear filtering allows one to reduce the effect of the noises
on the signals of the accelerometer $y_{a}$ and the linear velocity
$y_{v}$.

Define $p_{n}\overset{\Delta}{=}y_{v}-\hat{x}_{1}=x_{1}-\hat{x}_{1}$.
Then the two-steps $n$-th order observer is given by 
\begin{equation}
\begin{cases}
\dot{\hat{x}}{}_{2}^{\prime} & =-S\left(y_{g}\right)\hat{x}{}_{2}^{\prime}-\frac{\alpha_{1}}{g_{0}}p_{2},\\
\dot{p}_{i} & =-S\left(y_{g}\right)p_{i}+p_{i+1},\;(i=2,\cdots,n-1)\\
\dot{\hat{x}}_{1} & =-S\left(y_{g}\right)\hat{x}_{1}+y_{a}+\sum_{i=2}^{n}\alpha_{i}p_{i}-g_{0}\hat{x}{}_{2}^{\prime},\\
\dot{\hat{x}}_{2} & =-S\left(y_{g}-\gamma S\left(\hat{x}_{2}\right)\hat{x}'_{2}\right)\hat{x}_{2}.
\end{cases}\label{eq:observer-nth-order}
\end{equation}
Here, the gains $\alpha_{i}$, ($i=1,\ldots,n$) are positive and
chosen so that the polynomial $s^{n}+\alpha_{n}s^{n-1}+\alpha_{n-1}s^{n-2}+\alpha_{n-2}s^{n-3}+...+\alpha_{2}s+\alpha_{1}$
is Hurwitz. Moreover, $\hat{x}_{1}$ and $\hat{x}_{2}$ are estimations
of $x_{1}$ and $x_{2}$ respectively and $\hat{x}'_{2}$ is an intermediate
estimation of $x_{2}$. Using the estimation errors defined as $\tilde{x}_{1}=x_{1}-\hat{x}_{1}=p_{n}$,
$\tilde{x}_{2}^{\prime}=x_{2}-\hat{x}'_{2}=p_{1}$ and $\tilde{x}_{2}\overset{\Delta}{=}x_{2}-\hat{x}_{2}$,
we get the error dynamics as 
\begin{equation}
\begin{cases}
\dot{p}_{1} & =-S\left(\omega\right)p_{1}+\frac{\alpha_{1}}{g_{0}}p_{2},\\
\dot{p}_{i} & =-S\left(\omega\right)p_{i}+p_{i+1},\;\left(i=2,\cdots,n-1\right)\\
\dot{p}_{n} & =-S\left(\omega\right)p_{n}-g_{0}p_{1}-\sum_{i=2}^{n}\alpha_{i}p_{i}\\
\dot{\tilde{x}}_{2} & =-S\left(\omega\right)\tilde{x}_{2}+\gamma S^{2}\left(\hat{x}_{2}\right)\tilde{x}_{2}-\gamma S^{2}\left(\hat{x}_{2}\right)p_{1}.
\end{cases}\label{error_dynamics-nth-order}
\end{equation}
To run the error analysis, we set $z_{p_{1}}\overset{\Delta}{=}Rp_{1}$
, $z_{p_{i}}\overset{\Delta}{=}\frac{\alpha_{1}}{g_{0}}Rp_{i}$ ($i=2,\cdots,n$)
and $z_{2}\overset{\Delta}{=}R\tilde{x}_{2}$. Then one gets 
\begin{equation}
\begin{cases}
\dot{z}_{p_{i}} & =z_{p_{i+1}},\;(i=1,\cdots,n-1)\\
\dot{z}_{p_{n}} & =-\sum_{i=1}^{n}\alpha_{i}z_{p_{i}},\\
\dot{z}_{2} & =\gamma S^{2}\left(e_{z}-z_{2}\right)\left(z_{2}-z_{p_{1}}\right).
\end{cases}\label{eq:error_dynamics1-nth-order}
\end{equation}
This new dynamics is autonomous as for the first order case and defines
a time-invariant ordinary differential equation (ODE) on $\varUpsilon_{n}$.
Similarly to the previous state estimator, if one defines the state
$\xi_{n}\overset{\Delta}{=}\left(z_{p_{1}},\cdots,z_{p_{n}},z_{2}\right)\in\varUpsilon_{n}$,
one can write (\ref{eq:error_dynamics1-nth-order}) as $\dot{\xi}_{n}=F_{n}\left(\xi_{n}\right)$
where $F_{n}$ gathers the right-hand side of (\ref{eq:error_dynamics1-nth-order})
and defines a smooth vector field on $\varUpsilon_{n}$.

Note that the first $n$ lines of (\ref{eq:error_dynamics1-nth-order})
constitute a separate tilt estimator defined in $\mathbb{R}^{3n}$,
which is similar in the case $n=2$ to the one provided in \cite{Martin2016arxiv}.
We show hereinafter the convergence and the performances of this estimation
which are similar to the two-steps\textcolor{brown}{{} }first-order
tilt estimator. 
\begin{thm}
\label{th2} The time-invariant ODE defined by (\ref{eq:error_dynamics1-nth-order})
verifies the same statements as in Theorem\ref{th1} up to changing
$\varUpsilon_{1}$ and $\varUpsilon_{1}^{*}$ by $\varUpsilon_{n}$
and $\varUpsilon_{n}^{*}$, the first zero in the equilibrium points
now belonging to $\mathbb{R}^{3n}$ and $\alpha_{1}$ changed by $(\alpha_{1},\dots,\alpha_{n})$. 
\end{thm}

The proof of the theorem is given in Section~\ref{app:proof-th2}.

The higher order of the convergence of $\hat{x}_{2}^{\prime}$ allows
one to improve the robustness of the estimation to noise by increasing
the order $n$, as shown in simulations.
\begin{rem}
In the above construction, it is worth noticing that the gain-coefficients
$\alpha_{i}$, $1\leq i\leq n$ can be chosen time-varying and the
above construction remains unchanged till (\ref{eq:error_dynamics1-nth-order}).
One can, therefore, use the observer given in \cite{Chitour2002}
which improves the performances of the complementary filter regarding
a possible peaking phenomenon and noise.
\end{rem}

\subsection{One-step tilt observer designed in \texorpdfstring{{}$\mathbb{R}^{3}\times\mathbb{S}^{2}$}{R3xS2}
\label{subsec:RxS}}

In this section, we provide a tilt observer first given in \cite{Benallegue2017humanoids},
which is a slight improvement of the estimator proposed in \cite{Hua2016automatica}.
It is designed in $\mathbb{R}^{3}\times\mathbb{S}^{2}$ in one step
by using the available measures $y_{g}$, $y_{a}$ and $y_{v}$ and
is given by 
\begin{equation}
\begin{cases}
\dot{\hat{x}}_{1} & =-S(y_{g})\hat{x}_{1}-g_{0}\hat{x}_{2}+y_{a}+\alpha\tilde{x}_{1},\\
\dot{\hat{x}}_{2} & =-S(y_{g}+\gamma S(\hat{x}_{2})\tilde{x}_{1})\hat{x}_{2},
\end{cases}\label{eq:RxS_observer}
\end{equation}
where $\alpha$ and $\gamma$ are positive scalar gains which verify
the condition $\gamma g_{0}\leq\alpha^{2}$ and $\hat{x}_{1}$ and
$\hat{x}_{2}$ are the estimations of $x_{1}$ and $x_{2}$.

The observer proposed in \cite{Hua2016automatica} is recalled as
\begin{equation}
\begin{cases}
\dot{\hat{x}}_{1} & =-S(y_{g})\hat{x}_{1}-g_{0}\hat{x}_{2}+y_{a}+k_{1}^{v}\tilde{x}_{1}-k_{2}^{v}S^{2}\left(\hat{x}_{2}\right)\tilde{x}_{1},\\
\dot{\hat{x}}_{2} & =-S(y_{g}+k_{1}^{r}S(\hat{x}_{2})\tilde{x}_{1})\hat{x}_{2},
\end{cases}\label{eq:RxS_observer_Hua}
\end{equation}
where $k_{1}^{v}$, $k_{2}^{v}$ and $k_{1}^{r}$ are positive scalar
gains which verify the condition $k_{1}^{r}g_{0}\leq k_{1}^{v}k_{2}^{v}$.

We can notice that our proposed observer can be obtained from the
one of \cite{Hua2016automatica} by taking $k_{2}^{v}=0$, $k_{1}^{v}=\alpha$,
$k_{1}^{r}=\gamma$ and the condition on the gains becomes $k_{1}^{r}g_{0}\leq\left(k_{1}^{v}\right)^{2}$
instead of $k_{1}^{r}g_{0}\leq k_{1}^{v}k_{2}^{v}$.

Using the errors $\tilde{x}_{1}\overset{\Delta}{=}x_{1}-\hat{x}_{1}$
and $\tilde{x}_{2}\overset{\Delta}{=}x_{2}-\hat{x}_{2}$ as well as
$\omega=y_{g}$, a time-differentiation of these expressions provides
us with the following error dynamics: 
\begin{equation}
\begin{cases}
\dot{\tilde{x}}_{1} & =-S(\omega)\tilde{x}_{1}-\alpha\tilde{x}_{1}-g_{0}\tilde{x}_{2},\\
\dot{\tilde{x}}_{2} & =-S(\omega)\tilde{x}_{2}-\gamma S^{2}(\hat{x}_{2})\tilde{x}_{1},
\end{cases}\label{eq:RxS_error_dynamics}
\end{equation}

The error dynamics of the observer of \cite{Hua2016automatica} is
given by 
\begin{equation}
\begin{cases}
\dot{\tilde{x}}_{1} & =-S(\omega)\tilde{x}_{1}-\alpha\tilde{x}_{1}-g_{0}\tilde{x}_{2}+k_{2}^{v}S^{2}\left(\hat{x}_{2}\right)\tilde{x}_{1},\\
\dot{\tilde{x}}_{2} & =-S(\omega)\tilde{x}_{2}-\gamma S^{2}(\hat{x}_{2})\tilde{x}_{1},
\end{cases}\label{eq:RxS_error_dynamics_Hua}
\end{equation}

To run the analysis of errors, we define $z_{i}\overset{\Delta}{=}R\tilde{x}_{i}$.
We notice also that $R(\tilde{x}_{2}+\hat{x}_{2})=e_{z}$ which leads
to $R\hat{x}_{2}=e_{z}-z_{2}$ , we obtain this new error dynamics
of our proposed observer 
\begin{equation}
\begin{cases}
\dot{z}_{1} & =-\alpha z_{1}-g_{0}z_{2},\\
\dot{z}_{2} & =-\gamma S^{2}(e_{z}-z_{2})z_{1},
\end{cases}\label{eq:RxS_error_dynamics_2}
\end{equation}

We do the same for the observer of \cite{Hua2016automatica}, we get
the error dynamics as 
\begin{equation}
\begin{cases}
\dot{z}_{1} & =-\alpha z_{1}-g_{0}z_{2}+k_{2}^{v}S^{2}\left(e_{z}-z_{2}\right)z_{1},\\
\dot{z}_{2} & =-\gamma S^{2}\left(e_{z}-z_{2}\right)z_{1}.
\end{cases}\label{eq:RxS_error_dynamics_2_Hua}
\end{equation}

These new error dynamics are autonomous. In fact, if we define the
following state vector $\xi\overset{\Delta}{=}\left(z_{1},z_{2}\right)$
and the state space $\varUpsilon_{1}$, we can write these errors
dynamics as $\dot{\xi}=F\left(\xi\right)$ where $F$ defines smooth
vector fields on $\varUpsilon_{1}$.

Almost global asymptotic stability with respect to the origin $\left(0,0\right)$
is obtained in the case of the proposed observer with the gain condition
$\gamma g_{0}\leq\alpha^{2}$ and the proof is conducted in the same
way as in \cite{Hua2016automatica} by considering only one Lyapunov
function candidate given by 
\begin{align}
V & \overset{\Delta}{=}\frac{\left\Vert \alpha z_{1}+g_{0}z_{2}\right\Vert ^{2}}{2}+g_{0}^{2}\frac{\left\Vert z_{2}\right\Vert ^{2}}{2}.\label{eq:Lyap-1}
\end{align}

The time derivative of (\ref{eq:Lyap-1}) in view of (\ref{eq:RxS_error_dynamics_2})
yields 
\begin{align}
\dot{V} & =-\alpha\left(1-G_{0}\right)\left\Vert \alpha z_{1}+g_{0}z_{2}\right\Vert ^{2}\nonumber \\
 & +\alpha g_{0}^{2}G_{0}z_{2}^{T}S^{2}(e_{z})z_{2}\nonumber \\
 & -\alpha G_{0}\left(\left(\alpha z_{1}+g_{0}z_{2}\right)^{T}(e_{z}-z_{2})\right)^{2},\label{eq:d_Lyap-1}
\end{align}
where $G_{0}=\frac{\gamma g_{0}}{\alpha^{2}}\leq1$. This Lyapunov
function makes the convergence analysis much easier than that given
in \cite{Hua2016automatica}\emph{.}

\section{Attitude estimation observer\label{sec:Invariant}}

The measurements of the magnetometer provide the direction of the
magnetic field expressed in the local frame of the sensor. Usually
most of the measurement is constituted with the earth natural magnetic
field, which provides bi-dimensional data on the attitude of the sensor,
providing then enough inputs to reconstruct the full attitude and
having then some redundancy with the accelerometer for tilt estimation.
However, sometimes due to the proximity of sources of interference,
the magnetometer's measurements could lack the necessary reliability
to let it influence the critical tilt estimation, but remains the
best measurement available to reconstruct the orientation around the
vertical direction. In this case the solution is to use an estimation
allowing to tune the influence of the magnetometer on the tilt.

\subsection{Design of the attitude observer}

Let $\hat{R}\in SO(3)$ denote the estimate of $R$. The proposed
non-linear observer takes advantage of the estimator of $x_{2}$ designed
into $\mathbb{R}^{3}$ given by (\ref{eq:observer-nth-order}) and
the attitude estimator proposed by Mahony et. al \cite{mahony2008nonlinear},
and it is given by 
\begin{equation}
\begin{cases}
\dot{\hat{R}} & =\hat{R}S(y_{g}-\sigma),\\
\sigma & =\rho_{1}S(\hat{R}^{T}e_{z})\hat{x}'_{2}+\rho_{2}S(\hat{R}^{T}m)y_{m}\\
 & \;\qquad+\mu\hat{R}^{T}e_{z}\left(\hat{R}^{T}e_{z}\right)^{T}S\left(\hat{R}^{T}m\right)y_{m}.
\end{cases}\label{eq:redundent-observer}
\end{equation}
where $\rho_{1}$, $\rho_{2}$ and $\mu$ are positive scalar gains
and $\hat{x}'_{2}$ is given by the the first stage of any order of
the two-step tilt estimator from Section~\ref{sec:Tilt-estimation},
for example with (\ref{eq:observer-1st-order}).

In the case where $\rho_{2}=0$, we recover an estimate of the total
rotation with decoupled tilt in an essentially similar way as that
of \cite{Hua2016automatica} where the magnetometer has no influence
on the tilt. On the contrary, if $\mu=0$, the corresponding estimator
is closer to that of \cite{mahony2008nonlinear} and the estimator
fully uses the redundancy.

Let $\tilde{R}=R\hat{R}^{T}$ be the attitude estimation error. A
time-differentiation of the expression of (\ref{eq:redundent-observer})
and the use of equation (\ref{eq:error_dynamics-2-1st-order}) provides
us with the following error dynamics: 
\begin{alignat}{1}
\dot{\tilde{R}} & =\tilde{R}S\left(\tilde{\sigma}\right),
\end{alignat}
where $\tilde{\sigma}$ is given by 
\begin{alignat}{1}
\tilde{\sigma}=\begin{array}{l}
\left(I+\frac{\mu}{\rho_{2}}e_{z}e_{z}^{T}\right)\left(\rho_{1}S\left(e_{z}\right)\tilde{R}^{T}e_{z}+\rho_{2}S\left(m\right)\tilde{R}^{T}m\right)\\
-\rho_{1}S(e_{z})\tilde{R}^{T}z_{p_{1}}.
\end{array}
\end{alignat}
Using unit-quaternions instead of elements of $SO(3)$, one associates
$Q$ and $\hat{Q}$ with the rotations $R$ and $\hat{R}$ respectively,
and similarly the unit-quaternion error $\tilde{Q}=(\tilde{q}_{0},\tilde{q})=Q\odot\hat{Q}^{-1}$
with the attitude estimation error $\tilde{R}$. Here, $\tilde{q}_{0}\in\mathbb{R}$
and $\tilde{q}\in\mathbb{R}^{3}$ are the scalar and the vector components
of $\tilde{Q}$ respectively. We can, therefore, write 
\begin{alignat}{1}
\tilde{R} & =I+2\tilde{q}_{0}S\left(\tilde{q}\right)+2S^{2}\left(\tilde{q}\right),\nonumber \\
\rho_{1} & S\left(e_{z}\right)\tilde{R}^{T}e_{z}+\rho_{2}S\left(m\right)\tilde{R}^{T}m=-2\left(\tilde{q}_{0}I-S\left(\tilde{q}\right)\right)W_{\rho}\tilde{q},
\end{alignat}
with $W_{\rho}\overset{\Delta}{=}-\rho_{1}S^{2}\left(e_{z}\right)-\rho_{2}S^{2}\left(m\right)$
being a positive-definite symmetric matrix (\cite{A.Tayebi2013},
Lemma 2).

Set $\varpi\overset{\Delta}{=}\left(\tilde{q}_{0}I-S\left(\tilde{q}\right)\right)W_{\rho}\tilde{q}$.
The error dynamics written as a quaternion error dynamics is now given
by 
\begin{equation}
\begin{cases}
\dot{\tilde{q}}_{0} & \hspace{-3mm}=\tilde{q}^{T}\left(I+\frac{\mu}{\rho_{2}}e_{z}e_{z}^{T}\right)\varpi,\\
 & \hspace{-6mm}+\frac{1}{2}\rho_{1}\tilde{q}^{T}S(e_{z})\left(I-2\tilde{q}_{0}S\left(\tilde{q}\right)+2S^{2}\left(\tilde{q}\right)\right)z_{p_{1}},\\
\dot{\tilde{q}} & \hspace{-3mm}=-\left(\tilde{q}_{0}I+S\left(\tilde{q}\right)\right)\left(I+\frac{\mu}{\rho_{2}}e_{z}e_{z}^{T}\right)\varpi\\
 & \hspace{-6mm}-\frac{1}{2}\rho_{1}\left(\tilde{q}_{0}I+S\left(\tilde{q}\right)\right)S(e_{z})\left(I-2\tilde{q}_{0}S\left(\tilde{q}\right)+2S^{2}\left(\tilde{q}\right)\right)z_{p_{1}}.
\end{cases}\label{eq:error_dynamics-attitude-1}
\end{equation}
The above equation together with the first equation of (\ref{eq:error_dynamics-2-1st-order})
define a time-invariant ordinary differential equation (ODE) and,
by considering the state $\xi\overset{\Delta}{=}\left(z_{p_{1}},\tilde{Q}\right)$
and the state space $\varUpsilon\overset{\Delta}{=}\mathbb{R}^{3}\times\mathbb{S}^{3}$,
one can write (\ref{eq:error_dynamics-2-1st-order}) and (\ref{eq:error_dynamics-attitude-1})
as $\dot{\xi}=F\left(\xi\right)$ where $F$ gathers the right-hand
side of (\ref{eq:error_dynamics-attitude-1}) and defines a smooth
vector field on $\varUpsilon$. We analyze this dynamics in the next
section.

\subsection{Stability analysis}

Let us consider the following positive-definite differentiable function
\begin{eqnarray}
V & \overset{\Delta}{=} & \frac{\rho_{1}^{2}}{\alpha_{1}}\Vert z_{p_{1}}\Vert^{2}+2\tilde{q}^{T}W_{\rho}\tilde{q},\label{eq:Lyapunov-1-1}
\end{eqnarray}
which is clearly radially unbounded. 
\begin{thm}
\label{th5} \textup{\emph{The time-invariant ODE defined by (\ref{eq:error_dynamics-attitude-1})
verifies the following.}}
\end{thm}

\begin{enumerate}
\item \emph{Its equilibrium points are}
\[
\begin{array}{ccl}
\widetilde{\Omega}_{1} & = & \left\{ (0,(\pm1,0))\right\} ,\\
\widetilde{\Omega}_{2} & = & \left\{ (0,(0,\pm v_{j\rho})),\:j=1,2,3\right\} ,
\end{array}
\]
where $v_{j\rho}$ are unit eigenvectors of $W_{\rho}$ for $1\leq j\leq3$. 
\item \emph{All trajectories of (\ref{eq:error_dynamics-attitude-1}) converge
to one of the equilibrium points defined in item 1.} 
\item \emph{The set equilibrium $\widetilde{\Omega}_{1}$ which corresponds
to the equilibrium point $\left(z_{p_{1}}=0,\tilde{R}=I\right)$ is
asymptotically stable with a domain of attraction containing the domain
\begin{equation}
V_{c}\overset{\Delta}{=}\left\{ \xi=\left(z_{p_{1}},\tilde{Q}\right)\in\varUpsilon\mid V\left(\xi\right)<2\lambda_{min}(W_{\rho})\right\} .
\end{equation}
} 
\item \emph{The equilibria of the set $\widetilde{\Omega}_{2}$ are unstable
and the system is almost globally asymptotically stable with regard
to $\widetilde{\Omega}_{1}$.} 
\end{enumerate}
The proof is given in Section~\ref{app-pf-th5}.

\begin{rem}
The magnetic field measurements $y_{m}$ can also be filtered using
an additional unconstrained state on the unit sphere in the same way
as done for the tilt in order to improve robustness to noise. 
\end{rem}

\textcolor{red}{{} }
\begin{rem}
In the above, we have chosen, for the simplicity of the analysis,
to estimate the intermediate state $\hat{x}'_{2}$ with the two-steps
first order state observer given by (\ref{eq:observer-1st-order}).
One can also rely on the two-steps $n^{th}$ order state observer
given in (\ref{eq:error_dynamics1-nth-order}). For the corresponding
stability analysis, one replaces the $\Vert z_{p_{1}}\Vert^{2}$ term
in the Lyapunov function $V$ given in (\ref{eq:Lyapunov-1-1}) by
$\psi_{n}^{T}P_{\alpha}\psi_{n}$ given in (\ref{lyap-final}). 
\end{rem}

\section{Simulations\label{sec:Simulations}}

We show hereinafter results of the estimators in a simulated environment.

\subsection{Signal generation and initialization}

\begin{figure}
\begin{centering}
\includegraphics[width=0.8\columnwidth]{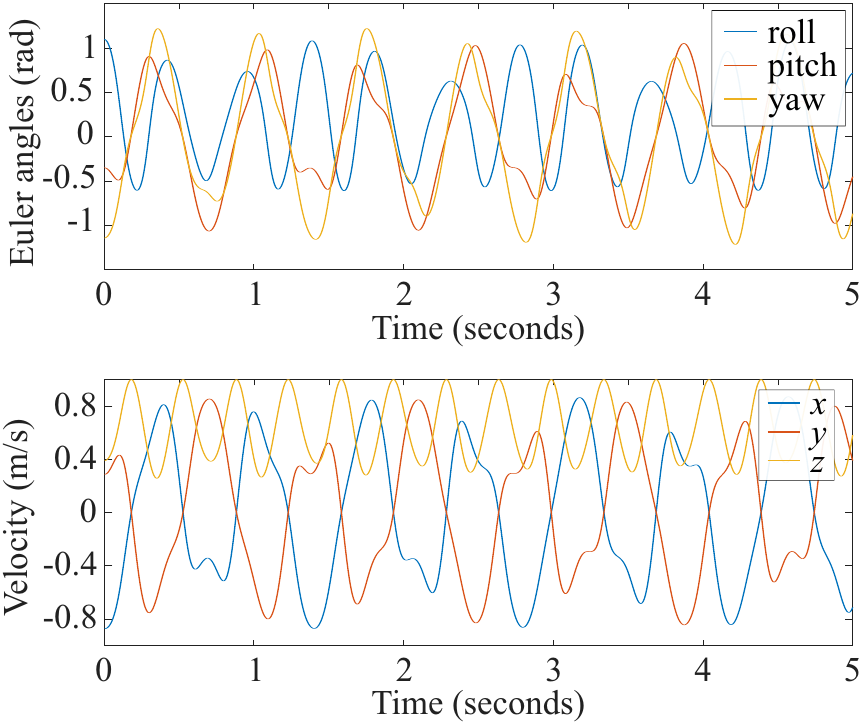}
\par\end{centering}
\caption{Plot showing the real state of the system. On the top the orientation
is shown in Euler angles and on the bottom, the velocity $x_{1}$
is shown in its three components.\label{fig:realstate}}
\end{figure}

In this section, we present simulation results showing the effectiveness
of the proposed estimators. We generated the signal $\omega$ and
$v$ with trigonometric functions and generated the trajectory of
$R$ by integration (see Figure \ref{fig:realstate}), then we simulated
the signals of the accelerometer $y_{a}$, the gyrometer $y_{g}$
and the magnetometer $y_{m}$ such that $m=\frac{1}{\sqrt{2}}(1,0,1)^{T}.$

These signals were used in two cases, ideal signals and noisy ones.
For the noisy signals, Gaussian noises were added to the four measurements,
the accelerometer $y_{a}$ the gyrometer $y_{g}$, the velocity sensor
$y_{v}$, and the normalized magnetometer $y_{m}$ to which a stronger
noise and a bias have been added to make it unreliable and unsuitable
to influence tilt estimation. The detail of the noise properties is
summarized in Table \ref{tab:noises}.

\begin{table}
\begin{centering}
\begin{tabular}{|c|c|c|}
\hline 
Measurement & Noise std. & Bias\tabularnewline
\hline 
Accelero. $y_{a}$ & 0.31 $\text{m/s}^{2}$ & $\left(0\;0\;0\right)^{T}$\tabularnewline
\hline 
Gyro. $y_{g}$ & 0.1 rad/s & $\left(0\;0\;0\right)^{T}$\tabularnewline
\hline 
Magneto. $y_{m}$ & 0.71 & $\left(0.2\;0.2\;0.2\right)^{T}$\tabularnewline
\hline 
Velocity $y_{v}$ & 0.31$\text{m/s}$ & $\left(0\;0\;0\right)^{T}$\tabularnewline
\hline 
\end{tabular}\smallskip{}
\par\end{centering}
\caption{Description of the noise parameters. \label{tab:noises}}
\end{table}

For each tested estimator the initial state was set to $\tilde{R}_{3}=2\left(\frac{m\times e_{z}}{\left\Vert m\times e_{z}\right\Vert }\right)\left(\frac{m\times e_{z}}{\left\Vert m\times e_{z}\right\Vert }\right)^{T}-I$
, which corresponds to an undesired equilibrium. The velocity estimation
was initialized to the current sensor value (for instance $\hat{x}_{1}(0)=x_{1}(0)$).

\subsection{Comparison between two stage tilt estimators}

The first test is to compare the tilt estimators presented in Section~\ref{sec:Tilt-estimation}.
Specifically, the first order, the second order and the third order
tilt estimators were compared for the perfect and the noisy measurements.
The estimators were designed to have the same (multiple) pole. The
parameters are detailed in Table~\ref{tab:Parameters-til}.
\begin{center}
\begin{table}
\begin{centering}
\begin{tabular}{|c|c|}
\hline 
Order & Parameters\tabularnewline
\hline 
\hline 
1st order & $\gamma=20$, $\alpha_{1}=2\sqrt{\gamma g_{0}}$\tabularnewline
\hline 
\multirow{1}{*}{2nd order} & $\gamma=20$, $\alpha_{1}=\left(2\sqrt{\gamma g_{0}}\right)^{2}$,
$\alpha_{2}=2\left(2\sqrt{\gamma g_{0}}\right)$\tabularnewline
\hline 
\multirow{2}{*}{3rd order} & $\gamma=20$, $\alpha_{1}=\left(2\sqrt{\gamma g_{0}}\right)^{3}$,\tabularnewline
 & $\alpha_{2}=3\left(2\sqrt{\gamma g_{0}}\right)^{2}$, $\alpha_{3}=3\left(2\sqrt{\gamma g_{0}}\right)$\tabularnewline
\hline 
\end{tabular}\smallskip{}
\par\end{centering}
\caption{Parameters of tested tilt estimators\label{tab:Parameters-til}}
\end{table}
\par\end{center}

\begin{figure}
\begin{centering}
\includegraphics[width=1\columnwidth]{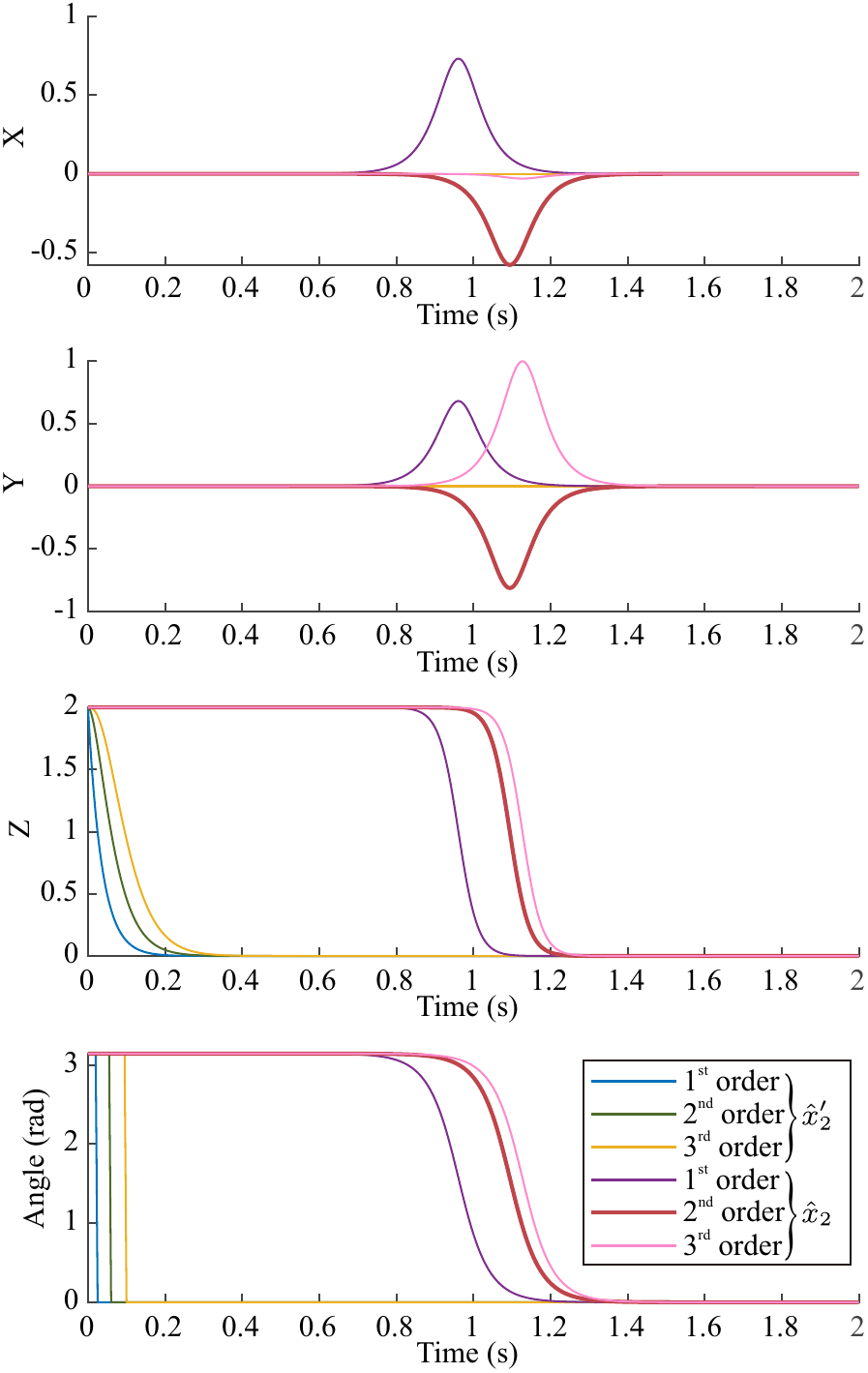}
\par\end{centering}
\caption{Plot showing the estimation error of the tilt vector $x_{2}$ for
the three orders of tilt estimators. For each order we show both the
intermediate estimation $\hat{x}'_{2}$ and the final one $\hat{x}{}_{2}$.
The three top images show the three components of the vector difference
error $z_{2}=R\left(x_{2}-\hat{x}_{2}\right)$ and the bottom plot
shows the evolution of the angle between the tilt $x_{2}$ and its
estimation $\hat{x}_{2}$. \label{fig:x_2-tilt}}
\end{figure}

The result of the simulation with perfect measurements is shown in
Figure~\ref{fig:x_2-tilt} where we compare the errors produced by
the estimations $\hat{x}_{2}$ but also the intermediate estimations
$\hat{x}'_{2}$. We see that the intermediate estimation errors converge
exponentially to zero while the estimation itself remains in the undesired
equilibrium. We can see that the first order estimator is obviously
the fastest followed by the other orders.

\begin{figure}
\begin{centering}
\includegraphics[width=1\columnwidth]{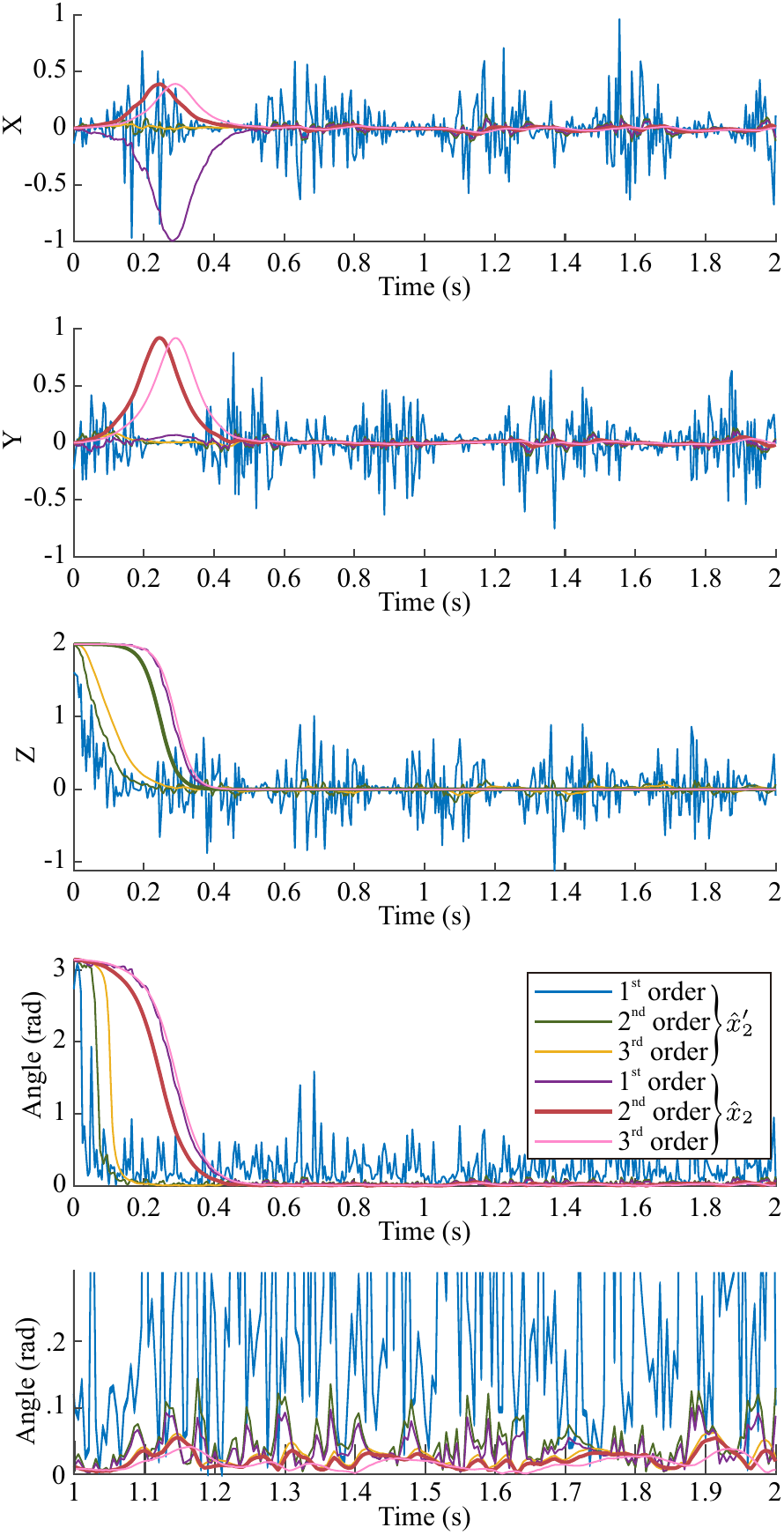}
\par\end{centering}
\caption{Plot showing the estimation error, in the case of noisy measurements,
of the tilt vector $x_{2}$ for the three orders of tilt estimators.
For each order we show both the intermediate estimation $\hat{x}'_{2}$
and the final one $\hat{x}{}_{2}$. The three top images show the
three components of the vector difference error $z_{2}=R\left(x_{2}-\hat{x}_{2}\right)$,
the 4th plot shows the evolution of the angle between the tilt $x_{2}$
and its estimation $\hat{x}_{2}$, and the bottom part shows an enlarged
plot of the second $[1,2]$ of the angle error.\label{fig:x_2-tilt_noisy}}
\end{figure}

However, the more interesting case of the noisy one displayed in Figure~(\ref{fig:x_2-tilt_noisy}).
We see then that with higher orders of the estimator better filtering
is provided. We see also that the sphere constraint of the final estimate
$\hat{x}_{2}$ allows to reduce the noise by removing the components
which are orthogonal to the constraints. Nevertheless the difference
between the second and the third order is small enough to consider
that the second order estimator is a good trade-off between complexity
and speed on one side and filtering quality on the other. Therefore,
in the following simulations we will use to feed the attitude estimator
in (\ref{eq:redundent-observer}) with $\hat{x}'_{2}$ and then compare
it with state-of-the-art approaches.

\subsection{Comparison between attitude estimators}

In this section, five estimators were compared.
\begin{enumerate}
\item The attitude estimator described in Section \ref{sec:Invariant},
with $\rho_{2}=0$ which means it is decoupled to avoid any impact
of the magnetometer on the tilt estimation. We refer to it by hierarchic.
\item The estimator described in Section \ref{sec:Invariant}, with $\mu=0$
using redundancy, referred to as Invariant.
\item The estimator in $\mathbb{R}^{3}\times\mathbb{S}^{2}$ in Section
\ref{subsec:RxS}, providing only tilt estimation, and referred to
as Benallegue 2017.
\item The estimator designed by Hua et al, named ``Observer 2'' in \cite{Hua2016automatica}
reported in (\ref{eq:RxS_observer_Hua}), that we refer to as Hua
2016.
\item The estimator described in the preprint \cite{Martin2016arxiv} by
Martin et al, named Martin 2016, which is based on an estimator equivalent
to 2nd order $\hat{x}'_{2}$ of Sec. \ref{subsec:State-observer-in-R3xR3xS2}
and another exponential estimator of the tilt and using TRIAD~\cite{shuster1981jgc}
to reconstruct the attitude. The estimation is designed for the tilt
to depend only on the accelerometer and the yaw angle only on the
magnetometer.
\end{enumerate}
Each estimator provides a specific tilt estimator. Note that the tilt
estimation of the invariant observer is the only one that requires
magnetometer's measurements.

The corresponding gains were designed to have the most equivalent
behavior possible, regarding their structure and the considered errors.
These estimators, as well as their tilt component and the gains used,
are summarized in Table~\ref{tab:estimators}\footnote{The tilt estimation column relates the different estimators to their
tilt estimation component and the gains column gives the gain values
adopting the notation used in each corresponding cited document.}.

\begin{table*}
\begin{centering}
{\small{}}%
\begin{tabular}{|c|c|c|}
\hline 
{\small{}Estimator} & {\small{}Tilt estimation} & {\small{}Gains}\tabularnewline
\hline 
\hline 
{\small{}Hierarchic (Observer (\ref{eq:redundent-observer}) with
$\rho_{2}=0$)} & {\small{}$\mathbb{R}^{3}\times\mathbb{R}^{3}\times\mathbb{S}^{2}$
(2nd order $\hat{x}{}_{2}$ of Sec. \ref{subsec:State-observer-in-R3xR3xS2})} & {\small{}$\gamma=20$, $\alpha=2\sqrt{\gamma g_{0}}=28.0143$,}\tabularnewline
\cline{1-2} \cline{2-2} 
{\small{}Invariant (Observer (\ref{eq:redundent-observer}) with $\mu=0$)} & {\small{}Invariant (uses magnetometer)} & {\small{}$\mu=20$ (hierarchic) or $\rho_{2}=20$ (invariant)}\tabularnewline
\cline{1-2} \cline{2-2} 
{\small{}Benallegue 2017 \cite{Benallegue2017humanoids} of Sec. \ref{subsec:RxS}} & {\small{}$\mathbb{R}^{3}\times\mathbb{S}^{2}$, (provides tilt only)} & {\small{}$\rho_{1}=20$, $\alpha_{1}=\gamma g_{0},$ $\alpha_{2}=2\sqrt{\gamma g_{0}}$}\tabularnewline
\hline 
{\small{}Hua 2016 \cite{Hua2016automatica}} & {\small{}Hua 2016 (Observer (\ref{eq:RxS_observer_Hua}))} & {\small{}$k_{1}^{v}=k_{2}^{v}=\alpha$, $k_{1}^{r}=k_{2}^{r}=\gamma$}\tabularnewline
\hline 
{\small{}Martin 2016 \cite{Martin2016arxiv}} & {\small{}$\mathbb{R}^{3}\times\mathbb{R}^{3}$ (2nd order $\hat{x}'_{2}$
of Sec. \ref{subsec:State-observer-in-R3xR3xS2})} & {\small{}$L=K=\frac{\alpha}{2}$, $M=\mu$}\tabularnewline
\hline 
\end{tabular}\smallskip{}
\par\end{centering}
\caption{Summary of compared estimators in the simulations.\label{tab:estimators}}
\end{table*}

\subsubsection{Perfect measurements}

\begin{figure}
\begin{centering}
\includegraphics[width=1\columnwidth]{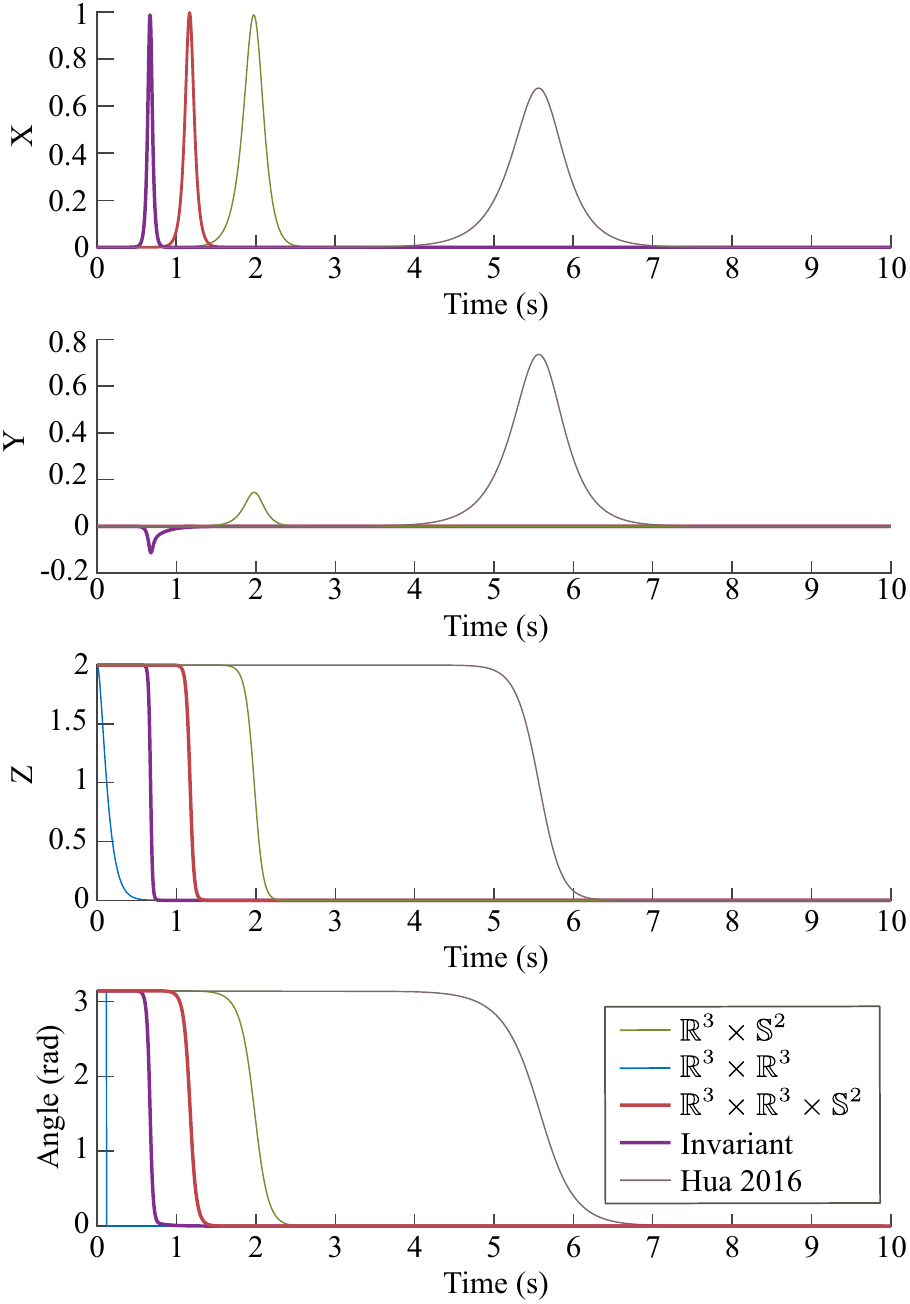}
\par\end{centering}
\caption{Plot showing the estimation error of the tilt vector $x_{2}$ for
the five tilt estimators. The three top images show the three components
of the vector difference error $z_{2}=R\left(x_{2}-\hat{x}_{2}\right)$
and the bottom plot shows the evolution of the angle between the tilt
$x_{2}$ and its estimation $\hat{x}_{2}$. Note that we use the names
of the second column of Table~\ref{tab:estimators}.\label{fig:x_2}}
\end{figure}

Figure \ref{fig:x_2} shows the evolution of the tilt error for the
five tilt estimators. The first estimator to converge is the one of
$\mathbb{R}^{3}\times\mathbb{R}^{3}$, namely 2nd order $\hat{x}'_{2}$
of Sec. \ref{subsec:State-observer-in-R3xR3xS2} which is not constrained
to the unit sphere, this is because the starting position is not an
equilibrium point for this vector. However, the normalization of this
vector gives a discontinuous trajectory visible at the bottom plot
showing the angle error. This is the estimation used in \cite{Martin2016arxiv}.
The next estimator to converge is the invariant one, this is due to
the fact that this estimator uses also the measurement of the magnetometer
to speedup the convergence. After that the estimator in $\mathbb{R}^{3}\times\mathbb{R}^{3}\times\mathbb{S}^{2}$
of Section \ref{subsec:State-observer-in-R3xR3xS2} is the next to
quickly converge, while staying continuous and constrained on the
unit sphere. The other estimators converge later, especially the estimation
of Hua 2016.

\begin{figure}
\begin{centering}
\includegraphics[width=1\columnwidth]{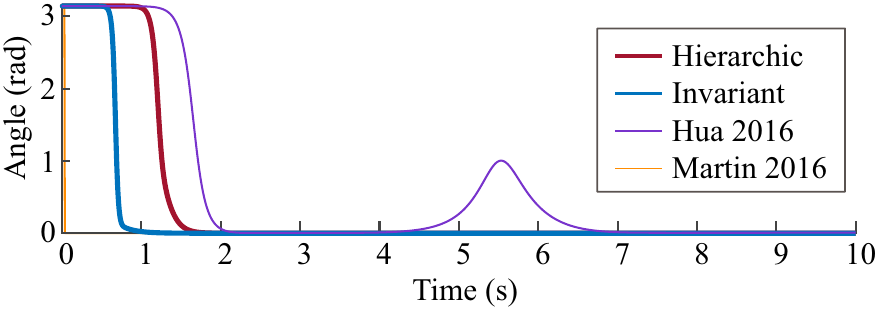}
\par\end{centering}
\caption{Plot showing the evolution of the angle between $R^{T}m_{p}$ and
its estimation $\hat{R}^{T}m_{p}$ for the five full attitude estimators.
\label{fig:ml}}
\end{figure}

Figure \ref{fig:ml} shows the evolution over time of the estimation
error of the vector $m_{p}\overset{\Delta}{=}e_{z}\times x_{3}\times e_{z}$
which is orthogonal to $e_{z}$ but pointing at the same horizontal
direction as $x_{3}$. The error is shown as an angle which can be
interpreted as a ``yaw angle error'' when the tilt error is small.
In this figure we see that the estimation of Martin 2016~\cite{Martin2016arxiv},
is discontinuous at another instant than the discontinuity of the
tilt, which means that the attitude had two discontinuities while
converging. The invariant estimator converges fast, taking full profit
from the redundancy. The Hierarchic, was the next estimator to converge.
We see finally that the estimation error of Hua 2016 moved at second
2 to zero. However, this does not correspond to the convergence of
the estimator since it took the tilt estimation 4 more seconds to
converge (see Figure \ref{fig:x_2}). This means that it only went
from an undesired equilibrium to another one. Note that some angles
increase and then decrease, and this happens because of the tilt estimation
converging at the same time and the orthogonality constraint being
respected.

\subsubsection{Noisy Measurements}

\begin{figure}
\begin{centering}
\includegraphics[width=1\columnwidth]{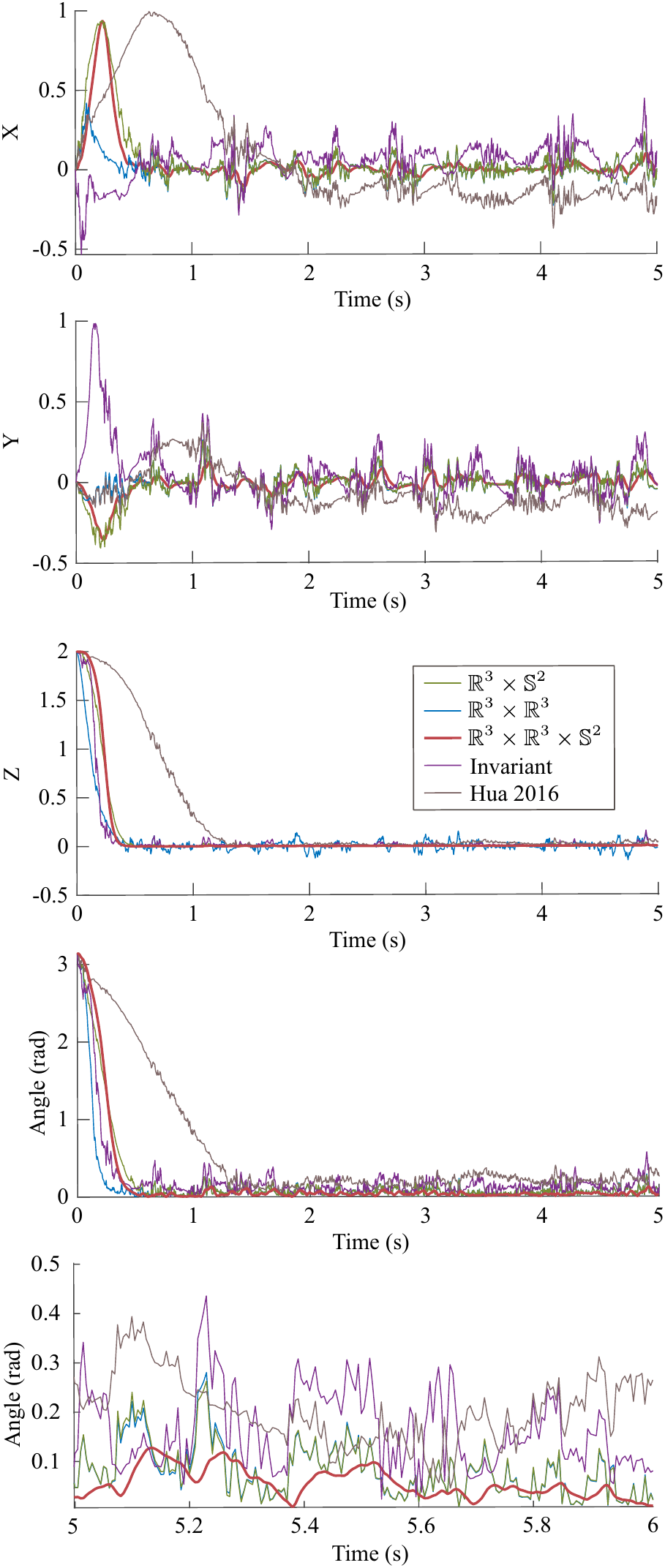}
\par\end{centering}
\caption{Plot showing the estimation error of the tilt vector $x_{2}$ for
the five tilt estimators while the measurements were noisy. The three
top images show the three components of the vector difference error
$z_{2}=R\left(x_{2}-\hat{x}_{2}\right)$, the 4th plot shows the evolution
of the angle between the tilt $x_{2}$ and its estimation $\hat{x}_{2}$,
and the bottom part shows an enlarged plot of the 6th second of the
angle error. Note that we use the names of the second column of Table~\ref{tab:estimators}.
\label{fig:x2_noise}}
\end{figure}

Figure \ref{fig:x2_noise} shows the tilt estimation error with the
difference and the angle, similarly to Figure \ref{fig:x_2}, with
an additional enlarged sample plot of the behavior after the convergence.

From this plot we see that the noise allowed the estimators to instantly
leave the repulsive undesired equilibrium. Then most estimators except
for Hua 2016 converge in less than half a second, the unconstrained
$\mathbb{R}^{3}\times\mathbb{R}^{3}$ being the fastest. After the
convergence of all the estimators, we see in the enlarged plot that
the estimator in $\mathbb{R}^{3}\times\mathbb{R}^{3}\times\mathbb{S}^{2}$
has the lowest tilt estimation error angle. The dynamics of the estimators
in $\mathbb{R}^{3}\times\mathbb{R}^{3}$ and $\mathbb{R}^{3}\times\mathbb{S}^{2}$
have identical local behavior near the desired equilibrium and are
almost superimposed in the steady behavior. Interestingly the invariant
observer gives worse estimations, that is because it involved the
unreliable magnetometer measurements which downgrade the performances.
We see in table \ref{tab:Average-tilt-error-noise} the mean value
of the tilt error angle over 8 seconds after the second 2 of the simulation.
The constrained $\mathbb{R}^{3}\times\mathbb{R}^{3}\times\mathbb{S}^{2}$
gives the best estimations and the invariant and Hua 2016 both give
the worst ones.

\begin{table}
\begin{centering}
\begin{tabular}{|c|c|}
\hline 
Tilt estimation & Mean tilt error angle\tabularnewline
\hline 
$\mathbb{R}^{3}\times\mathbb{R}^{3}\times\mathbb{S}^{2}$ & 0.0442 rad\tabularnewline
\hline 
Invariant & 0.1543 rad\tabularnewline
\hline 
$\mathbb{R}^{3}\times\mathbb{S}^{2}$ & 0.0748 rad\tabularnewline
\hline 
Hua 2016 & 0.1960 rad\tabularnewline
\hline 
$\mathbb{R}^{3}\times\mathbb{R}^{3}$ & 0.0749 rad\tabularnewline
\hline 
\end{tabular}\smallskip{}
\par\end{centering}
\caption{Average tilt error angles during 8 seconds after the convergence of
the estimators.\label{tab:Average-tilt-error-noise}}
\end{table}

\begin{table}
\begin{centering}
\begin{tabular}{|c|c|}
\hline 
Attitude est. & Mean $m_{p}$ error angle\tabularnewline
\hline 
Hierarchic & 0.2374 rad\tabularnewline
\hline 
Invariant & 0.2511 rad\tabularnewline
\hline 
Hua 2016 & 0.2671 rad\tabularnewline
\hline 
Martin 2016 & 0.3036 rad\tabularnewline
\hline 
\end{tabular}\smallskip{}
\par\end{centering}
\caption{Average $R^{T}m_{p}$ estimation error angles during 8 seconds after
the convergence of the estimators.\label{tab:Average-e_x-error-noise}}
\end{table}

\begin{figure}
\begin{centering}
\includegraphics[width=1\columnwidth]{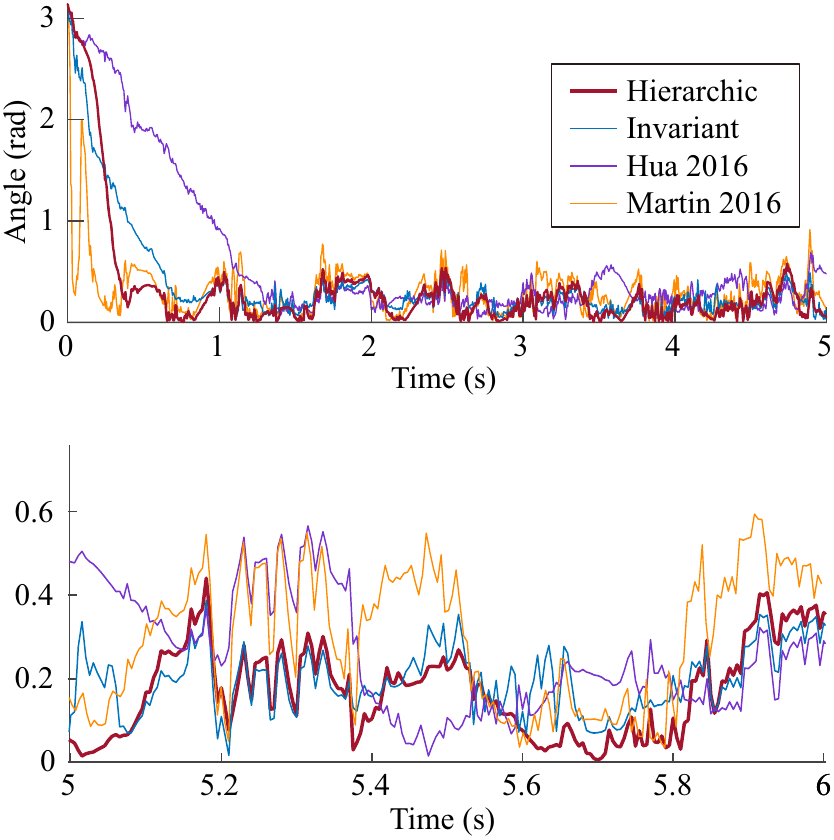}
\par\end{centering}
\caption{Plot showing the evolution of the angle between $R^{T}m_{p}$ and
its estimation $\hat{R}^{T}m_{p}$ for the five full attitude estimators
under noisy measurements. \label{fig:m_noise}}
\end{figure}

Figure \ref{fig:m_noise} shows the evolution of the estimation of
$R^{T}m_{p}$ together with a zoom on the 6-th second of the simulation.
We see that the estimations converge in the first second except for
Hua 2016. The high level of noise in the magnetometer produces a poor
estimation quality, but in the steady behavior a difference can be
shown between observers. This can be quantitatively assessed by looking
at Table \ref{tab:Average-e_x-error-noise}\textbf{ }showing the average
error angle values in the interval {[}2s,10s{]}\textbf{. }Martin 2016
has low quality estimations because the estimation of yaw is performed
independently from the measurements of the accelerometer. The other
estimators take profit from the better reliability of the tilt estimation
and provide a relatively similar level of performance with a slight
advantage to the hierarchic estimator.

Note that a behavior between the hierarchic and the invariant estimator
can be obtained by choosing values of $\mu$ and $\rho_{2}$ appropriately,
especially that small values of $\rho_{2}$ provide better theoretical
convergence guarantees without downgrading excessively the quality
of the estimation.

\section{Conclusion\label{sec:Discussion-and-conclusion}}

We have presented a set of attitude estimators using the measurements
of an accelerometer, a gyrometer, a magnetometer, and a linear velocity
expressed in the local frame. These estimators are intended to be
used in specific cases, mostly related to the reliability of the magnetometer
for tilt estimation. Indeed, the magnetometer can be either reliable,
unreliable, or totally unavailable. For instance, an invariant complementary
filter has good performances when the magnetometer is reliable but
is disturbed when it is not. Among the estimators, we developed a
second order complementary filter for the tilt, and we augmented it
with an attitude estimator allowing to tune how much we use the magnetometer
measurement in the tilt estimation. We have assessed the performances
of these estimators through simulations of perfect and noisy measurements.

\section{Appendix\label{sec:appendix}}

\subsection{Proof of Theorem~\ref{th1}}

\label{app:proof-th1} 1) One easily checks that the time-invariant
ODE defined by (\ref{eq:error_dynamics-2-1st-order}) leaves invariant
$\varUpsilon_{1}$ since along its trajectories, $e_{z}-z_{2}$ keeps
a constant norm equal to one. Moreover, it admits two equilibrium
points $(0,0)$ and $(0,2e_{z})$. \\
 {} Let us now consider the following positive-definite differentiable
function $V_{1}:\varUpsilon_{1}\rightarrow\mathbb{R}^{+}$ 
\begin{eqnarray}
V_{1} & \overset{\Delta}{=} & \frac{1}{2\alpha_{1}}z_{p_{1}}^{T}z_{p_{1}}+\frac{1}{2\gamma}z_{2}^{T}z_{2}\label{eq:Lyapunov-1st-order}
\end{eqnarray}
then the time derivative of $V_{1}$ is given by 
\begin{align*}
\dot{V}_{1} & =-\Vert z_{p_{1}}\Vert^{2}+z_{2}^{T}S^{2}\left(e_{z}-z_{2}\right)z_{2}-z_{2}^{T}S^{2}\left(e_{z}-z_{2}\right)z_{p_{1}}
\end{align*}
If we use $w=S(e_{z})z_{2}$, we can write 
\[
\dot{V}_{1}\leq-\left[\begin{array}{cc}
\left\Vert z_{p_{1}}\right\Vert  & \Vert w\Vert\end{array}\right]\left[\begin{array}{cc}
1 & -\frac{1}{2}\\
-\frac{1}{2} & 1
\end{array}\right]\left[\begin{array}{c}
\left\Vert z_{p_{1}}\right\Vert \\
\Vert w\Vert
\end{array}\right]\leq0,
\]
and $\dot{V}_{1}<0$ if $(z_{p_{1}},w)\neq0$. This is equivalent
to $\left(z_{p_{1}},z_{2}\right)$ not being an equilibrium point,
i.e., all trajectories of\emph{ }(\ref{eq:error_dynamics-2-1st-order})
converge to one of the two equilibrium points defined previously.
\\
 {} 2) The linearized system at $(0,0)$ is Hurwitz yielding that
$(0,0)$ is locally exponentially stable. At $(0,2e_{z})$, the linearized
system clearly admits two positive real eigenvalues. We can conclude
that the system (\ref{eq:error_dynamics-2-1st-order}) is almost globally
asymptotically stable with respect to the origin ($0,0$). Moreover,
the set of points of $\varUpsilon_{1}$ for which $V_{1}$ has values
less than $V_{1}(0,2e_{z})=2/\gamma$ is clearly included in the basin
of attraction of ($0,0$). \\
 {} 3) Let $K$ be a compact set in $\varUpsilon_{1}^{*}$ and $\varrho>0$.
It is easy to see that for every $(z_{p_{1}},z_{2})\in K$ one has
that 
\[
V_{1}(z_{p_{1}},z_{2})\leq\frac{C_{1}(K)}{2\alpha_{1}}+\frac{2}{\gamma}C_{2}(K),
\]
with $C_{2}(K)<1$. By the remark at the end of the argument of item
2) and by choosing $\alpha_{1}$ large enough, one gets $K$ is in
the basin of attraction of ($0,0$). \\
 {} To obtain the last statement, we first prove that there exists
$C_{3}(K)>0$ such that $\Vert w\Vert\leq C_{3}(K)\Vert z_{2}\Vert$
for $(z_{p_{1}},z_{2})\in K$. It is enough to check that for $z_{2}$
small enough. If one writes $z_{2}=(e_{z}^{T}z_{2})e_{z}+z_{2}^{\perp}$
one gets that $\Vert z_{2}\Vert\leq2\Vert z_{2}^{\perp}\Vert$ (by
using that $z_{2}\in\mathbb{S}_{e_{z}}$) and $w=z_{2}^{\perp}$ for
$z_{2}$ small, hence the claim. Next, one deduces that there exists
$C(K)>0$ such that 
\[
\dot{V}_{1}\leq-C(K)(\Vert z_{p_{1}}\Vert^{2}+\Vert z_{2}\Vert^{2}),
\]
for trajectories starting in $K$ (and staying in a compact neighborhood
of $K$ in the basin of attraction of ($0,0$)). Then $\dot{V}_{1}\leq-2C(K)\min(\alpha_{1},\gamma)V_{1}$.
By taking $\alpha_{1}$ and $\gamma$ large enough, one gets the conclusion. 

\subsection{Proof of Theorem~\ref{th2}}

\label{app:proof-th2} The argument is similar to that of Theorem~\ref{th1}.
For that purpose consider the Hurwitz $n\times n$ matrix in companion
form $A_{\alpha}=J_{n}-ae_{n}^{T}$, $J_{n}$ stands for the $n$-th
Jordan block, $a=(\alpha_{1},\dots,\alpha_{n})^{T}$ and $e_{n}=(0,\cdots,0,1)^{T}$.
Then set $M_{\alpha}=A_{\alpha}\otimes I_{3\times3}$ and $\psi_{n}=\left(z_{p_{1}},\dots,z_{p_{n}}\right)\in\mathbb{R}^{3n}$.
Note that the $n$ first equations in (\ref{eq:error_dynamics1-nth-order})
can be written $\dot{\psi}_{n}=M_{\alpha}\psi_{n}$. Let $P_{\alpha}$
be the positive definite real symmetric matrix, unique solution of
the Lyapunov equation 
\[
M_{\alpha}^{T}P_{\alpha}+P_{\alpha}M_{\alpha}=-I_{3n\times3n}.
\]
Recall that $\Vert P_{\alpha}\Vert\leq\frac{C}{re_{\alpha}}$, where
$C_{n}$ is a universal positive constant and $re_{\alpha}>0$ is
the minimum of $-Re(\lambda)$, $Re$ stands for the real part and
$\lambda$ is any eigenvalue of the $A_{\alpha}$, cf. \cite{horn_johnson_1991}.
\\
 One now considers the Lyapunov function 
\begin{equation}
V_{n}=\psi_{n}^{T}P_{\alpha}\psi_{n}+\frac{1}{2\gamma}z_{2}^{T}z_{2}.\label{lyap-final}
\end{equation}
We now follow exactly the argument of Theorem~\ref{th1} and replace
$\Vert z_{p_{1}}\Vert$ by $\Vert\psi_{n}\Vert$ to get the conclusion. 

\subsection{Proof of Theorem~\ref{th5}}

\label{app-pf-th5} Let us prove the four items of the theorem.

1) The equilibria are calculated by solving the equation $\dot{\xi}=0$.
The solutions of this equation system are given by $\left(z_{p_{1}}=0,\varpi=0\right)$.
We know from (\cite{A.Tayebi2013}, Lemma 3) that $\varpi=0$ is equivalent
to $(\tilde{q}_{0}=\pm1,\tilde{q}=0)$ or $(\tilde{q}_{0}=0,\tilde{q}=\pm v_{\rho})$
where $v_{\rho}$ is one of the unit eigenvectors of $W_{\rho}$.
This completes the proof of item 1.

2) Using the error dynamics given by (\ref{eq:error_dynamics-attitude-1}),
the time derivative of $V$ is then given by 
\begin{equation}
\dot{V}=-2\rho_{1}^{2}\Vert z_{p_{1}}\Vert^{2}+4\tilde{q}^{T}W_{\rho}\dot{\tilde{q}},
\end{equation}
which can be developed into 
\begin{alignat}{1}
\dot{V}= & -2\rho_{1}^{2}\Vert z_{p_{1}}\Vert^{2}\\
 & -4\tilde{q}^{T}W_{\rho}\left(\tilde{q}_{0}I+S\left(\tilde{q}\right)\right)\left(I+\frac{\mu}{\rho_{2}}e_{z}e_{z}^{T}\right)\varpi\nonumber \\
 & +2\rho_{1}\tilde{q}^{T}W_{\rho}\left(\tilde{q}_{0}I+S\left(\tilde{q}\right)\right)S(e_{z})\tilde{R}^{T}z_{p_{1}}.
\end{alignat}
Using the definition of the vector $\varpi$, we obtain 
\begin{align*}
\dot{V} & =-2\rho_{1}^{2}\Vert z_{p_{1}}\Vert^{2}-4\Vert\varpi\Vert^{2}-4\frac{\mu}{\rho_{2}}(e_{z}^{T}\varpi)^{2}\\
 & \qquad+2\rho_{1}\varpi^{T}S(e_{z})\tilde{R}^{T}z_{p_{1}}
\end{align*}
which can be bounded with the following expression 
\begin{align}
\dot{V} & \leq-2\rho_{1}^{2}\left\Vert z_{p_{1}}\right\Vert ^{2}-4\left\Vert \varpi\right\Vert ^{2}+2\rho_{1}\left\Vert \varpi\right\Vert \left\Vert z_{p_{1}}\right\Vert .
\end{align}

The right-hand side of the above inequality is a quadratic form in
$(\Vert z_{p_{1}}\Vert,\Vert\varpi\Vert)$ which is clearly negative
definite. One easily verifies that $\dot{V}<0$ if $\xi=\left(z_{p_{1}},\tilde{Q}\right)$
is not an equilibrium. Since (\ref{eq:error_dynamics-attitude-1})
is autonomous and $V$ is radially unbounded, one can use Lasalle's
invariance theorem. Therefore, every trajectory converges asymptotically
to a trajectory along which $\dot{V}\equiv0$.

Since $V$ is non-increasing, $V\left(\xi\right)<2\lambda_{min}(W_{\rho})$
at $t=0$, implies that $\left\Vert \tilde{q}(t)\right\Vert <1$ for
every $t\geq0$. Since the trajectory converges to one of the equilibrium
points, it must be one with ($z_{p_{1}}=0,\tilde{q}=0$) which corresponds
to $\widetilde{\Omega}_{1}$ because this is the only one contained
in $V_{c}$.

4) The undesired equilibria characterized by $\tilde{q}_{0}=0$ are
given by $X=(0,(0,v_{\rho}))$. Let us show that $X=(z_{p_{1}}=0,(\tilde{q}_{0}=0,\tilde{q}=v_{\rho}))$
is unstable. The linearized error dynamics around the unstable equilibrium
$X=(0,(0,v_{\rho}))$ is given by 
\begin{equation}
\dot{\xi}=A\xi,\label{eq:A}
\end{equation}
with $A$ given in Equation (\ref{Avalue}) in the box \vpageref{Avalue}.

\begin{algorithm*}[!t]
Let us consider $v_{\rho}^{\bot}=S\left(v_{\rho}\right)e_{z}$, we
can write the matrix $A$ as
\begin{equation}
A=\left[\begin{array}{ccccc}
-\alpha_{1}I &  & 0 &  & 0\\
\\
-\frac{1}{2}\rho_{1}v_{\rho}^{\bot T} &  & \lambda_{\rho}\left(1+\frac{\mu}{\rho_{2}}\left(v_{\rho}^{T}e_{z}\right)^{2}\right) &  & -\frac{\mu}{\rho_{2}}\left(v_{\rho}^{T}e_{z}\right)\left(v_{\rho}^{\bot}\right)^{T}\left(\lambda_{\rho}I-W_{\rho}\right)\\
\\
\frac{1}{2}\rho_{1}S\left(v_{\rho}\right)S\left(e_{z}\right)\left(I-2v_{\rho}v_{\rho}^{T}\right) &  & -\lambda_{\rho}\frac{\mu}{\rho_{2}}\left(e_{z}^{T}v_{\rho}\right)v_{\rho}^{\bot} &  & \left(I+\frac{\mu}{\rho_{2}}v_{\rho}^{\bot}v_{\rho}^{\bot T}\right)\left(\lambda_{\rho}I-W_{\rho}\right)
\end{array}\right]\label{Avalue}
\end{equation}
\end{algorithm*}

It is clear that there is at least one positive eigenvalue of the
matrix $A$. Thus, there exists an unstable manifold of dimension
at least one in neighborhoods of the $\widetilde{\Omega}_{2}=\left\{ (0,(0,\pm v_{j\rho})),\:j=1,2,3\right\} $,
and since all trajectories converge to an equilibrium point, then
(\ref{eq:error_dynamics-attitude-1}) is almost globally asymptotically
stable with respect to the two equilibrium points $\widetilde{\Omega}_{1}=\left\{ (0,(\pm1,0))\right\} $
which correspond to $\left(z_{p_{1}}=0,\tilde{R}=I\right)$. This
completes the proof.

\bibliographystyle{plain}
\bibliography{biblio-local}

\begin{thebibliography}{10}

\bibitem{allibert2014estimating}
Guillaume Allibert, Dinuka Abeywardena, Moses Bangura, and Robert Mahony.
\newblock Estimating body-fixed frame velocity and attitude from inertial
  measurements for a quadrotor vehicle.
\newblock In {\em 2014 IEEE Conference on Control Applications (CCA)}, pages
  978--983. IEEE, 2014.

\bibitem{allibert2016velocity}
Guillaume Allibert, Robert Mahony, and Moses Bangura.
\newblock Velocity aided attitude estimation for aerial robotic vehicles using
  latent rotation scaling.
\newblock In {\em 2016 IEEE International Conference on Robotics and Automation
  (ICRA)}, pages 1538--1543. IEEE, 2016.

\bibitem{Benallegue2017humanoids}
Mehdi Benallegue, Abdelaziz Benallegue, and Yacine Chitour.
\newblock {Tilt estimator for 3D non-rigid pendulum based on a tri-axial
  accelerometer and gyrometer}.
\newblock In {\em 2017 IEEE-RAS 17th International Conference on Humanoid
  Robotics (Humanoids)}, pages 830--835. IEEE, nov 2017.

\bibitem{Bloesch-RSS-12}
Michael Bloesch, Marco Hutter, Mark Hoepflinger, Stefan Leutenegger, Christian
  Gehring, C~David Remy, and Roland Siegwart.
\newblock {State Estimation for Legged Robots - Consistent Fusion of Leg
  Kinematics and {\{}IMU{\}}}.
\newblock In {\em Proceedings of Robotics: Science and Systems}, Sydney,
  Australia, jul 2012.

\bibitem{Chitour2002}
Y.~Chitour.
\newblock Time-varying high-gain observers for numerical differentiation.
\newblock {\em IEEE Trans. on Automatic Control}, 47(9), September 2002.

\bibitem{franklin1994feedback}
Gene~F Franklin, J~David Powell, and Abbas Emami-Naeini.
\newblock {\em {Feedback control of dynamic systems}}, volume~3.

\bibitem{horn_johnson_1991}
Roger~A. Horn and Charles~R. Johnson.
\newblock {\em Topics in Matrix Analysis}.
\newblock Cambridge University Press, 1991.

\bibitem{Hua2010cep}
Minh-Duc Hua.
\newblock {Attitude estimation for accelerated vehicles using GPS/INS
  measurements}.
\newblock {\em Control Engineering Practice}, 18(7):723--732, jul 2010.

\bibitem{Hua2017cdc}
Minh-Duc Hua, Tarek Hamel, and Claude Samson.
\newblock {Riccati nonlinear observer for velocity-aided attitude estimation of
  accelerated vehicles using coupled velocity measurements}.
\newblock In {\em 2017 IEEE 56th Annual Conference on Decision and Control
  (CDC)}, pages 2428--2433. IEEE, dec 2017.

\bibitem{Hua2016automatica}
Minh-Duc Hua, Philippe Martin, and Tarek Hamel.
\newblock {Stability analysis of velocity-aided attitude observers for
  accelerated vehicles}.
\newblock {\em Automatica}, 63:11--15, jan 2016.

\bibitem{mahony2008nonlinear}
Robert Mahony, Tarek Hamel, and Jean-Michel Pflimlin.
\newblock Nonlinear complementary filters on the special orthogonal group.
\newblock {\em IEEE Transactions on automatic control}, 53(5):1203--1217, 2008.

\bibitem{Martin2008ifac}
Philippe Martin and Erwan Sala{\"{u}}n.
\newblock {An Invariant Observer for Earth-Velocity-Aided Attitude Heading
  Reference Systems}.
\newblock In {\em IFAC Proceedings Volumes}, volume~41, pages 9857--9864.
  Elsevier, jan 2008.

\bibitem{martin2008invariant}
Philippe Martin and Erwan Sala{\"u}n.
\newblock An invariant observer for earth-velocity-aided attitude heading
  reference systems.
\newblock {\em IFAC Proceedings Volumes}, 41(2):9857--9864, 2008.

\bibitem{Martin2010ICRA}
Philippe Martin and Erwan Salaun.
\newblock {The true role of accelerometer feedback in quadrotor control}.
\newblock In {\em 2010 IEEE International Conference on Robotics and
  Automation}, pages 1623--1629. IEEE, may 2010.

\bibitem{Martin2016cdc}
Philippe Martin and Ioannis Sarras.
\newblock {A semi-global model-based state observer for the quadrotor using
  only inertial measurements}.
\newblock In {\em 2016 IEEE 55th Conference on Decision and Control (CDC)},
  pages 7123--7128. IEEE, dec 2016.

\bibitem{Martin2016arxiv}
Philippe Martin, Ioannis Sarras, Minh-Duc Hua, and Tarek Hamel.
\newblock {A global exponential observer for velocity-aided attitude
  estimation}.
\newblock {\em arXiv preprint arXiv:1608.07450}, aug 2016.

\bibitem{mifsud:hal-01142399}
Alexis Mifsud, Mehdi Benallegue, and Florent Lamiraux.
\newblock {Estimation of Contact Forces and Floating Base Kinematics of a
  Humanoid Robot Using Only Inertial Measurement Units}.
\newblock In {\em IEEE/RSJ International Conference on Intelligent Robots and
  Systems (IROS 2015)}, page 6p., Hamburg, Germany, sep 2015.

\bibitem{Roberts2011cdc}
Andrew Roberts and Abdelhamid Tayebi.
\newblock {On the attitude estimation of accelerating rigid-bodies using GPS
  and IMU measurements}.
\newblock In {\em IEEE Conference on Decision and Control and European Control
  Conference}, pages 8088--8093. IEEE, dec 2011.

\bibitem{shuster1981jgc}
Malcolm~David Shuster and S~D\_ Oh.
\newblock Three-axis attitude determination from vector observations.
\newblock {\em Journal of guidance and Control}, 4(1):70--77, 1981.

\bibitem{A.Tayebi2013}
A.~Tayebi, A.~Roberts, and A.~Benallegue.
\newblock Inertial vector measurements based velocity-free attitude
  stabilization.
\newblock {\em IEEE Transactions on Automatic Control}, 58(11):2893--2898, Nov
  2013.

\bibitem{Wieber2016}
Pierre-Brice Wieber, Russ Tedrake, and Scott Kuindersma.
\newblock {\em {Modeling and Control of Legged Robots}}, pages 1203--1234.
\newblock Springer International Publishing, Cham, 2016.

\end{thebibliography}

\end{document}